 \documentclass{article}

\usepackage{amssymb, amsmath}
\usepackage{latexsym}
\usepackage{graphics}
\usepackage{url}
\usepackage{enumerate}
\usepackage{amsfonts}
\usepackage{geometry}
\usepackage{graphicx}
\usepackage{amscd}
\usepackage{amssymb}
\usepackage{amsmath}
\usepackage[latin1]{inputenc}
\usepackage{epsfig,pstricks,pst-node,pst-tree,ifthen}
\usepackage{amsfonts}
\usepackage[all]{xy}
\usepackage{tikz}
\usetikzlibrary{fit}
\usepackage[all]{xy}
\usepackage{graphicx}
\usepackage{float}
\usepackage{amsrefs}
\usepackage{amsthm}
\usepackage{nicematrix}
\usepackage{esstixfrak}

\newtheorem{theo}{Theorem}
\newtheorem{coro}[theo]{Corollary}

\newtheorem{lemm}[theo]{Lemma}
\newtheorem{prop}[theo]{Proposition}

\newenvironment{proo}[1][\proofname]{\normalfont{\itshape
#1{:}}\quad\mdseries\ignorespaces}
{{$\Box$}{\vskip\belowdisplayskip}}
\newtheorem{defi}[theo]{Definition}
\newtheorem{rem}[theo]{Remark}

\newcommand{\N}{{\mathbb N}}

\newcommand{\K}{{\mathbb K}}

\begin{document}

\title{Commutators and commutator subgroups of the Riordan group}

\author{Ana Luz\'on, Manuel A. Mor\'on and Luis Felipe Prieto-Mart\'inez}

\maketitle

\begin{abstract}
We calculate the derived series of the Riordan group. To do that, we study a nested sequence of its subgroups, herein denoted by $\mathcal G_k$. By means of this sequence, we first obtain the $n$-th commutator subgroup of the Associated subgroup. This fact allows us to get some related results about certain groups of formal power series and to complete the proof of our main goal, Theorem \ref{main.2} in this paper.

\end{abstract}


\section{Introduction}







If $\mathcal G$ is a group and $g, h\in \mathcal G$, $[g,h]=g^{-1}h^{-1}gh $ is the commutator of $g$ and $h$. Let $\mathcal C=\{[g,h], | \, g,h\in \mathcal G\}$ be the subset of all commutators of $\mathcal G$. Denote, as usual, by $[\mathcal G, \mathcal G]$ to the commutator subgroup of $\mathcal G$, that is, $[\mathcal G, \mathcal G]$ is the subgroup generated by the set $\mathcal C$. We can define iteratively the $n$-th commutator subgroup, denoted by $\mathcal G^{(n)}$, as  $\mathcal G^{(1)}=[\mathcal G, \mathcal G]$ and $\mathcal G^{(n)}=[\mathcal G^{(n-1)}, \mathcal G^{(n-1)}]$. The sequence of groups $\mathcal G^{(n)}$ is called \textit{the derived series of $\mathcal G$}. Recall that $[\mathcal G,\mathcal G]$ is the smallest normal subgroup $\mathcal N$ such that the quotient group $\mathcal G/\mathcal N$ is abelian.

Following Guralnick \cite{G}, the problem of determining when $[\mathcal G, \mathcal G]=\mathcal C$ is of particular interest.  There are some examples in the literature where the commutator subgroup is not exclusively formed by commutators, see \cite{C, I, G} and the references therein. Ore in \cite{O.C} studied the commutator subgroup of certain groups. He also conjectured in \cite{O.C}
that any element in a finite non-abelian simple group is a commutator. This conjecture has been stablished in 2010 in \cite{LIEBECK}. In particular, for these kind of groups, all elements in the commutator subgroup are commutators. Thompson dealed with the corresponding problem for the special and general linear groups of finite matrices in \cite{Th}. Recently, it has been studied for certain groups of infinite matrices \cite{GO, S}.

Related to all above, the main result in this paper is:

 \begin{theo} \label{main.2} Let $\mathbb K$ be a field of characteristic 0. For $n\geq 1$:
\begin{equation}\label{equation.main.2}\mathcal R^{(n)}=\{(d,h)\in\mathcal R:d\in(1+x^{2^n-n}\mathbb K[[x]]), h\in (x+x^{2^n}\mathbb K[[x]])\}\end{equation}

\noindent and all of its elements are commutators of elements in $\mathcal R^{(n-1)}$.

\end{theo}

\noindent $\mathcal R$ above and in the rest of this paper denotes the Riordan group with entries in a  field $\mathbb K$ of characteristic 0 and $\mathcal R^{(n)}$, for $n\geq 1$, is the $n$-th commutator subgroup of $\mathcal R$. So, we describe the whole derived series of the Riordan group. A problem related to Theorem \ref{main.2} was posed by Shapiro in \cite{Sh.OQ}. In particular, our Theorem \ref{main.2} for $n=2$ gives the whole answer to Shapiro's problem.

As a consequence of the above statement, the matrices in $\mathcal R^{(n)}$ have all the elements in the main diagonal equal to 1 and the following $2^n-n-1$ diagonals  have all their entries equal to 0. That is, they are of the form
$$\begin{bNiceMatrix} 1 \\ 
 \Ddots[line-style=solid] &\Ddots \\ 
0 & \\ 
\Vdots & \Ddots &  \\
 0& & & &1  \\ 
\Ddots[line-style=solid] & \Ddots & & 0&  \\
 \times &  & & & \Ddots  & & \Ddots \\ \times  & \times &   &0 &  & 0&& 1  \\  \Vdots &\Vdots  & \Ddots  &\Ddots  & \Ddots & & \Ddots  & & \Ddots\\ \times  & \times  &    &  \times  &   & 0 & & & \\ \Vdots& \Vdots  & \Ddots  &   &\Ddots  &  &\Ddots  &  &    \end{bNiceMatrix} $$


The Riordan group $\mathcal R$  is a group whose elements are a special type of infinite lower triangular matrices with entries in $\mathbb K$, called Riordan matrices.   The  Riordan group was introduced, under this name and in a more restrictive context, by L. Shapiro and collaborators in  \cite{SHAPIRO.ORIGINAL}. Since then, many authors have developed this topic \cite{H, LM.U, MRSV,  S.CS, WW, Z}. Some previous related results can be found in \cite{BBN, JAB1, JAB2, ROGERS, S.C, VS}. The concept was initially motivated by its applications to Combinatorics \cite{S.CS}. During the last few years, the Riordan group and some related topics have been studied in \cite{BH, CLW, C.G1, CLMPS.L, HHS, JN, LMP.GGI, LMR, P, Y}. For the purpose of this paper, it is of special interest the work \cite{LMMPS.IL}, where the Riordan group  is described as an inverse limit of an inverse sequence of groups of finite matrices.







The derived series is also an important aspect to understand the algebraic structure of the Riordan group. Moreover, analogous description of commutators subgroups for the  Riordan groups of finite matrices introduced in \cite{LMMPS.IL} is also provided as a consequence of Theorem \ref{main.2}. These groups of finite matrices turn out to be solvable and, consequently, the Riordan group is a pro-solvable group.

In this work, two subgroups of the Riordan group, see Section \ref{sect.basics}, play an important role:  the Toeplitz or Appell subgroup $\mathcal T$ and the associated or Lagrange subgroup $\mathcal A$. They are also of interest because of their relationships with certain groups of formal power series.


Denote by $\mathbb K[[x]]$ to the set of formal power series $d=d_0+d_1x+d_2x^2+\ldots $ with coefficients in $\mathbb K$. As it is detailed in the following section, each Riordan matrix can be identified with a pair of formal power series $(d,h)$. As a consequence, the Riordan group $\mathcal R$ and its subgroups are also related to some groups of elements in $\mathbb K[[x]]$.  Consider $\mathcal F_0\subset \mathbb K[[x]]$ to be the multiplicative group of units of formal power series. That is, $\mathcal F_0$ is the multiplicative group of formal power series of order 0: $f_0+f_1x+f_2x^2+\ldots \text{ with }f_0\neq 0 $. Consider also $\mathcal F_1\subset \mathbb K[[x]]$ to be the group of invertible formal power series with respect to the composition. In this case, $\mathcal F_1$ is the  group of all series of order 1: $g_1x+g_2x^2+\ldots$ ($g_1\neq 0 $ and $g_0=0$) with the composition.  In this work, we denote by $g^{-1}$ to the compositional inverse of $g\in\mathcal F_1$. Finally, denote by $\mathcal J$   to   the subgroup of $\mathcal F_1$ of formal power series of the form: $x+g_2x^2+\ldots $ which is known in the bibliography as the \emph{substitution group} of formal power series (see \cite{JENN}). We refer to \cite{BB.S} for an interesting survey concerning to the algebraic structure of $\mathcal F_1$ and $\mathcal J$. In this paper, using the matricial approach, we obtain and improve some of the results in \cite{BB.S} related to the derived subgroup of $\mathcal F_1$. The groups $\mathcal F_0,\mathcal F_1,\mathcal J$ are isomorphic to certain subgroups of $\mathcal R$ (as explained in Remark \ref{rem.equiv}). Let us note that $\mathcal R$ is isomorphic to a semidirect product $\mathcal F_0 \rtimes \mathcal F_1 $.

Let $\mathbb F_p$ be the finite field consisting of $p$ elements for $p$ prime.  If we  consider the set of formal power series $\mathbb F_p[[x]]$  instead of $\mathbb K[[x]]$, the object analogue to the substitution group has also been extensively studied in the literature and it is known as the  \emph{Nottingham group} (see \cite{C.NG} for an overview of this topic).



The structure of this article is the following. We review the relevant basic facts about the Riordan group  in Section \ref{sect.basics}. In Section \ref{sect.GK},  we introduce a family of nested subgroups of $\mathcal R$, denoted by $\mathcal G_k$, that helps us to understand the proof and statement of the main theorem. We also think that this sequence of subgroups is interesting itself. It is described in Theorem \ref{prop.gknested}. We also prove therein a result, Proposition \ref{FE.MCWS1},  concerning a functional equation in the context of formal power series (Weighted Schr\"oder Equation). This result, which is itself interesting, is both necessary for  Theorem \ref{main.2} and a good example to understand the technique used in the rest of the proofs of this work: \emph{induction in groups of finite matrices}. In Section \ref{sec.fps}, we prove that $\mathcal A^{(n)}=\mathcal G_{2^n}$ (Theorem \ref{main.1bis}) and, as consequence, we give a description of the $n$-th commutator subgroup of $\mathcal F_1$ (Theorem \ref{main.1}). In Section \ref{sec.WSE}, we combine Proposition \ref{FE.MCWS1} and Theorem \ref{main.1} to prove Theorem \ref{main.2}. Finally, we discuss the consequences of theorems \ref{main.2} and \ref{main.1} for the  groups of finite Riordan matrices and we include some comments concerning future work.

\section{Basic facts about formal power series and Riordan matrices} \label{sect.basics}

In this paper $\mathbb N$ represents the set $\{0,1,2,3,\cdots\}\subset\mathbb K$,
$[x^k]$ denotes the $k$-th coefficient in the series expansion.

The concept of a  {\it Riordan matrix} and the related definition of  {\it Riordan group} appeared in the foundational paper \cite{SHAPIRO.ORIGINAL} due to Shapiro, Getu, Woan and Woodson. The original definition of a Riordan matrix given in \cite{SHAPIRO.ORIGINAL} is more restrictive than that used currently in the literature, which is precisely that we are going to use herein. A Riordan matrix is a matrix $D=(d_{i,j})_{i,j\in\N}$   whose columns are the coefficients of successive terms of a
geometric progression, in $\K[[x]]$, where the initial term $d$  is a
formal power series of order $0$ and with common ratio $h$, where $h$ is a formal
power series of order $1$.  From now on, we denote such a matrix $D$  by $(d,h)$. An element $(d,h)=(d_{i,j})_{i,j\in\mathbb N}$ in the Riordan group $\mathcal R$
  is an infinite matrix whose
entries are $d_{i,j}=[x^i]d(x)h^j(x)$.
Note that, by definition, these matrices are invertible
infinite lower triangular. The set of all matrices $(d,h)$ with the usual product of matrices forms a group called the \emph{Riordan group} which is denoted  by $\mathcal R$. In terms of the involved formal power series, the operations in the group are:
\[
(d,h)(l,m)=(d l(h),m(h)),  \qquad
(d,h)^{-1}=\left(\frac{1}{d(h^{-1})},h^{-1}\right)
\]

\noindent For any Riordan matrix $(d,h)$ there is a formal power series $A=\sum_{i=0}^\infty a_{i}x^i$ of order 0, called the A-sequence of $(d,h)$, with the property
\begin{equation} \label{eq.patternAred} d_{i,j}=\sum_{k=0}^{i-j}a_kd_{i-1,j-1+k} \qquad i,j\geq1 \end{equation}
It is known that $h^{-1}(x)=\frac{x}{A(x)}$ or, equivalently, $h(x)=x \cdot A(h(x))$. See \cite{L, ROGERS}. See also page 401 in \cite{LMP.I} for a proof of the existence and features of the A-sequence depending only on general group theoretic properties. The action induced by $(d,h)$ in $\mathbb K[[x]]$ is given by
\[
(d,h)\alpha=d\alpha(h) \quad \text{for}\ \alpha\in\mathbb K[[x]].
\]

\noindent  Many authors call the above equality the Fundamental Theorem of Riordan Matrices. The expression $(d,h)\alpha$ writen matricially corresponds to:
$$\begin{bmatrix} d_{00}\\ d_{10} & d_{11} \\ d_{20}& d_{21} & d_{22} \\ \vdots & \vdots  & \vdots  \ddots \end{bmatrix}\begin{bmatrix} \alpha_0\\ \alpha_1 \\ \alpha_2 \\ \vdots \end{bmatrix} \text{ where }\alpha=\alpha_0+\alpha_1x+\alpha_2x^2+\ldots$$

\noindent  Note that $(d,h)$ is a weighted composition operator in $\mathbb K[[x]]$.

For every $n\in\mathbb N$ consider the general linear group $GL(n+1,\mathbb K)$
formed by all $(n+1)\times(n+1)$ invertible matrices with
coefficients in $\mathbb K$. Since every Riordan matrix is lower triangular we have a natural
homomorphism $\Pi_n: \mathcal{R}\rightarrow GL(n+1,\mathbb K)$ given by $
\Pi_n((d_{i,j})_{i,j\in\mathbb N})=(d_{i,j})_{i,j=0,1,\cdots,n} $, as considered  in \cite{LMMPS.IL}. We denote by $(d,h)_n=\Pi_n((d,h))$ and by $\mathcal R_n=\Pi_n(\mathcal R)$. In the sequel we refer to these groups as Riordan groups of finite matrices.


We
can recover the group $\mathcal{R}$ as the inverse limit of the
inverse sequence of groups
$\{(\mathcal{R}_n)_{n\in\mathbb N},(P_n)_{n\in\mathbb N}\}$ where
$P_n:\mathcal R_{n+1}\rightarrow\mathcal R_n$ is such that if $D\in\mathcal R_{n+1}$,
$P_n(D)$ is obtained from $D$ by deleting its last row and its last
column, i.e.
$P_{n}((d_{i,j})_{i,j=0,1,\cdots,n+1})=(d_{i,j})_{i,j=0,1,\cdots,n}$.
See again \cite{LMMPS.IL}. Obviously if $n=0$ then $\mathcal R_0=\mathbb K^*$ with the
usual product in $\mathbb K^{*}=\mathbb K\setminus\left\{0\right\}$.

We now describe the two subgroups announced in the introduction.  First,  we have the \emph{Toeplitz subgroup} $\mathcal T$, made up by the elements of the type $(d,x)$. This subgroup is normal and abelian. Second, consider also the \emph{associated subgroup}  $\mathcal A$, formed by the elements of the type $(1,h)$. Recall that $\mathcal R$ is isomorphic to a semidirect product $\mathcal T \rtimes \mathcal A $. Another important fact relating the Riordan group with some groups of formal power series is:

\begin{rem} \label{rem.equiv} The groups $\mathcal F_0$ and $\mathcal F_1$  described in the introduction (and so any of their subgroups) are isomorphic to some subgroups of the Riordan group.
In fact, there is a natural isomorphism between $\mathcal F_0$ and the \emph{Toeplitz subgroup} $\mathcal T$ given by
$$f\longmapsto (f,x) $$
\noindent and another one between $\mathcal F_1$ and the \emph{associated subgroup} $\mathcal A$
$$g\longmapsto (1,g^{-1}). $$
\end{rem}

\section{The groups $\mathcal G_k$} \label{sect.GK}

This section is devoted to understand some aspects about a family of nested subgroups of $\mathcal A$, with a band of null diagonals  under the main one, that are present in the proofs of the main theorems of this article. These groups are also related to the balls with respect to the ultrametric introduced in \cite{LM.U}.  Moreover, the group $\mathcal G_{k+1}$ is isomorphic, via the isomorphism described in Remark \ref{rem.equiv}, to the groups $\mathfrak{G}_{k}$  introduced by Jennings in \cite{JENN}.

\begin{defi}  For $k\geq 2$, we define $\mathcal G_k=\{(1,h): h\in (x+ x^k\mathbb K[[x]])\}$.

\end{defi}

The following proposition clarifies the shape of the matrices in each $\mathcal G_k$:

\begin{theo}  \label{prop.gknested} $ $

\begin{itemize}

\item[(1)]  $\{\mathcal G_k\}_{k\geq 2}$ is a nested family of normal subgroups of $\mathcal G_2$.

\item[(2)] Let $(1,h)=(d_{ij})_{0\leq i,j<\infty}\in \mathcal G_k$, such that $h=x+h_kx^k+h_{k+1}x^{k+1}+\ldots$ Then 
$$d_{j+m,j}=\begin{cases}1&\text{if } m=0\\ 0 &\text{if } 1\leq m\leq k-2\\ j h_{j+1} &\text{if } k-1\leq m\leq 2k-2 \end{cases} $$

 \item[(3)]  Matrices in $\mathcal G_k$ have an A-sequence of the type $(1,0,\ldots,0,\alpha_{k-1},\alpha_{k},\ldots)$ where $\alpha_{l-1}= h_{l}$  for $1\leq l\leq 2k-2$.

\item[(4)] Let $(d,h)\in\mathcal R$, such that $(1,h)\in\mathcal G_k$ and $d=d_0+d_1x+\ldots$ Then:
$$d_{j+m,j}=\begin{cases}d_0&\text{if } m=0\\ d_m &\text{if } 1\leq m\leq k-2\\ d_m+j h_{j+1} &\text{if } k-1\leq m\leq 2k-2 \end{cases} $$

\end{itemize}

\end{theo}

\begin{proo} Consider the homomorphisms $\left.\Pi_l\right|_{\mathcal G_2}:\mathcal G_2\to \mathcal R_l$ defined in Section \ref{sect.basics}. Then $\mathcal G_k=ker(\left.\Pi_{k-1}\right|_{\mathcal G_2})$, and so all of them are normal subgroups of $\mathcal G_2$. The subgroups are nested because the diagram

\[\xymatrix{
\mathcal G_2 \ar[d]^{\Pi_{k-1}} \ar[drr]^{\Pi_{k}}\\
 \mathcal R_{k-1}  & & \mathcal R_{k} \ar[ll]^{P_{k-1}}
} \]
\noindent is commutative.

Note that conditions in (2) give rise to the following picture of the matrix $(1,h)$:
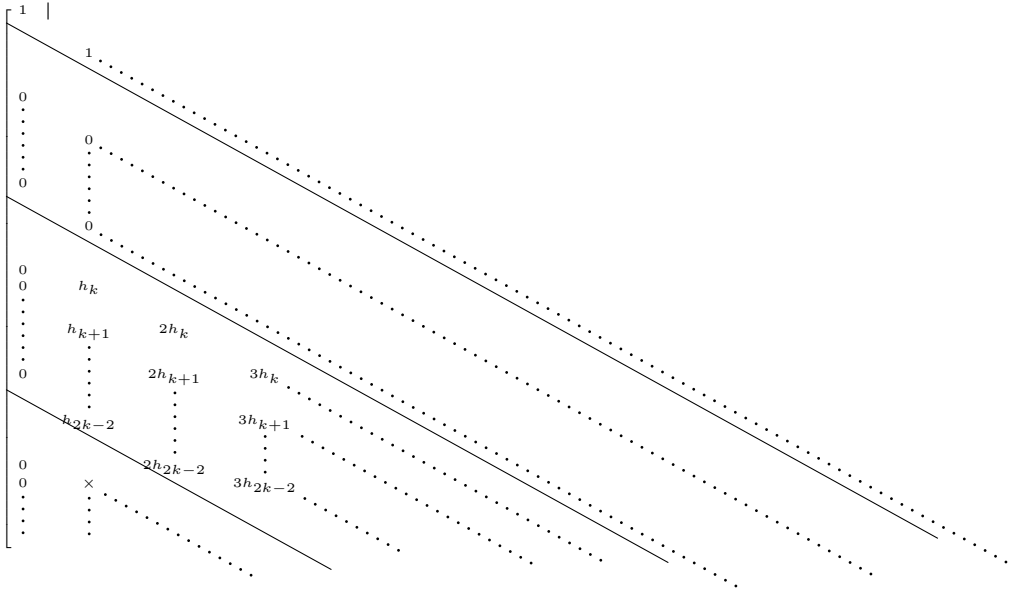
\begin{tiny}
$$\begin{bNiceArray}{c | c c c c c c c c c c c c c c } 1 & \textcolor{white}{aaaa} &\textcolor{white}{aaaa} &\textcolor{white}{aaaa} &\textcolor{white}{aaaa} &\textcolor{white}{aaaa} &\textcolor{white}{aaaa} &\textcolor{white}{aaaa} &\textcolor{white}{aaaa} &\textcolor{white}{aaaa} &\textcolor{white}{aaaa} &\textcolor{white}{aaaa} &\textcolor{white}{aaaa} &\textcolor{white}{aaaa} &\textcolor{white}{aaaa}\\ 
 \Ddots[line-style=solid] & 1\\ 
0 & & \Ddots \\ 
\Vdots & 0 & &   \\
 0&  \Vdots  &\Ddots  & & \\ 
\Ddots[line-style=solid] & 0 & & & & \\
 0 &  & \Ddots & & \Ddots  & &  \\
 0  & h_{k} &   &&  & && \\
 \Vdots & h_{k+1} & 2h_{k} &  & & &  & &  \\
0& \Vdots & 2h_{k+1} & 3h_{k}\\
\Ddots[line-style=solid] & h_{2k-2}  & \Vdots &  3h_{k+1} & \Ddots \\
0&  &2h_{2k-2}  & \Vdots &\Ddots & & & & & \\
0 & \times &&3h_{2k-2} & &&&&\\
\Vdots & \Vdots& \Ddots & &\Ddots & & & & & & & & & \\
\textcolor{white}{a}&\textcolor{white}{a} &\textcolor{white}{a} &\textcolor{white}{a} &\textcolor{white}{a} &\textcolor{white}{a} &\textcolor{white}{a} &\textcolor{white}{a} &\textcolor{white}{a} &\textcolor{white}{a} &\textcolor{white}{a} &\textcolor{white}{a} &\textcolor{white}{a} &\textcolor{white}{a} &\textcolor{white}{a}  \end{bNiceArray} $$
\end{tiny}

\noindent  The main diagonal has all its elements equal to 1. After this, the matrix has a band of $k-2$ null diagonals, followed by  a band of $k-1$ diagonals which entries form an arithmetic progression. This is a consequence of being $(1,h)$ a Riordan matrix.

Concerning the third part of the result, since the generating function of the A-sequence is  $\frac{x}{h^{-1}}$, and $h\in (x+x^k\K[[x]])$ the statement is clear.






Note that the A-sequence of $(d,h)$ is the same as that of $(1,h)$. Then, by the construction pattern of Riordan matrices in terms of the A-sequence, we can finally prove (4) as a consequence of (2).

\hspace{15cm}\end{proo}

For $i\geq 1$, let us denote by $H_i$ the map $H_i:\mathcal R\to\mathbb K$, given by $H_i((d,h))=[x^i]h$.  First we have:

\begin{lemm} \label{lemm.H} For $k\geq 2$ and $2\leq i\leq 2k-2$,  the restriction $\left.H_i\right|_{\mathcal G_k}$ is a homomorphism between $ \mathcal G_k$  and the abelian group $(\mathbb K,+)$.
\end{lemm}

\begin{proo}  Let $(1,v),(1,w)\in \mathcal G_2$, where $v(x)=x+v_2x^2+v_3x^3+\ldots$ and $w(x)=x+w_2x^2+w_3x^3+\ldots$ The product $(1,v)\cdot (1,w)$ is equal to:
\begin{equation} \label{eq.unadiagmas} \begin{bmatrix}1 \\ 0 &1\\ 0 & v_2 & 1 \\ 0 & v_3 & 2 v_2 & 1 \\  \vdots & \vdots & \vdots & \vdots & \ddots   \end{bmatrix} \begin{bmatrix}1 \\ 0 &1\\ 0 & w_2 & 1 \\ 0 & w_3 & 2 w_2 & 1 \\  \vdots & \vdots & \vdots & \vdots & \ddots   \end{bmatrix}= \begin{bmatrix}1 \\ 0 &1\\ 0 & v_2+w_2 & 1 \\ 0 & v_3+w_3+2v_2w_2 & 2 (v_2+w_3) & 1 \\  \vdots & \vdots & \vdots & \vdots & \ddots   \end{bmatrix} \end{equation}

\noindent From this expression we can easily check that $H_2((1,v)\cdot (1,w))=H_2((1,v))+H_2((1,w))$.

Let us now consider the general case  $k\geq 3$. In order to see that  for $2\leq i\leq 2k-2$, the restriction $\left.H_i\right|_{\mathcal G_k}$ is an homomorphism, we distinguish two cases.

For the case $2\leq i\leq k-1$, according to Theorem \ref{prop.gknested}, the elements in $\mathcal G_k$ have a band of $k-2$ null diagonals after the main one.  This already shows that,  $\left.H_i\right|_{\mathcal G_k}$ is a constant map (the trivial homomorphism).

For the case $k\leq i\leq 2k-2$, take the multiplication of the following two matrices in $\mathcal G_k$:
$$(1,v)(1,w)=\begin{bmatrix}1 \\ 0 & 1 \\ 0 & 0 & 1 \\ \vdots & \vdots &  \ddots& \ddots \\ 0 & 0 & \hdots & 0 & 1\\ 0 & v_k & 0 & \hdots & 0& 1 \\ 0 & v_{k+1} & 2v_k & 0 &  \hdots & 0& 1 \\ \vdots & \vdots & \vdots & \ddots & \ddots & &  \ddots \end{bmatrix}\begin{bmatrix}1 \\ 0 & 1 \\ 0 & 0 & 1 \\ \vdots & \vdots &  \ddots& \ddots \\ 0 & 0 & \hdots & 0 & 1\\ 0 & w_k & 0 & \hdots & 0& 1 \\ 0 & w_{k+1} & 2w_k & 0 &  \hdots & 0& 1 \\ \vdots & \vdots & \vdots & \ddots & \ddots & &  \ddots \end{bmatrix}= $$
\begin{equation}\label{eq.demcomm}=\begin{bmatrix}1 \\ 0 & 1 \\ 0 & 0 & 1 \\ \vdots & \vdots &  \ddots& \ddots \\ 0 & 0 & \hdots & 0 & 1\\ 0 & v_k+w_k & 0 & \hdots & 0& 1 \\ 0 & v_{k+1}+w_{k+1} & 2(v_k+w_k) & 0 &  \hdots & 0& 1 \\ \vdots & \vdots & \vdots & \ddots & \ddots & &  \ddots \end{bmatrix} \end{equation}

 \noindent  Denote $(1,v)= (a_{lj})_{0\leq l,j<\infty}$, $(1,w)= (b_{lj})_{0\leq l,j<\infty}$ and $(1,v)\cdot (1,w)= (c_{lj})_{0\leq l,j<\infty}$.    By the rule of multiplication of matrices, since $(1,w)\in \mathcal G_k$, $b_{i1}=w_i$ and using Theorem \ref{prop.gknested} we have:
$$c_{i1}= H_i((1,v)) (1,w))= \sum_{j=1}^ia_{ij}b_{j1}=a_{i1}+\sum_{j=k}^ia_{ij}w_j =v_i+w_i+\sum_{j=k}^{i-1} (i-j)v_{i-j+1}w_j $$


\noindent Finally, since $(1,v)\in \mathcal G_k$,  $c_{i1}=v_i+w_i$.



\hspace{15cm}\end{proo}

At this moment, we want to point out that  Lemma \ref{lemm.H} implies the commutativity of the subgroups $\Pi_{2k-1}(\mathcal G_{k})<\mathcal R_{2k-1}$ and so  we can get  that $[\mathcal G_k,\mathcal G_k]<\mathcal G_{2k-1}$. More accurately, we have:

\begin{theo}\label{theo.partialcommute} {For any two matrices $(1,v),(1,w)\in \mathcal G_k$, the finite Riordan matrices $(1,v)_{2k-1}$ and $(1,w)_{2k-1}$ commute.}
\end{theo}

\begin{proo}  Denote $(1,v)= (a_{lj})_{0\leq l,j<\infty}$, $(1,w)= (b_{lj})_{0\leq l,j<\infty}$,  $(1,v) (1,w)= (c_{lj})_{0\leq l,j<\infty}$,  $(1,w) (1,v)= (d_{lj})_{0\leq l,j<\infty}$,  where $v(x)=x+v_kx^k+\ldots$ and $w(x)=x+w_kx^k+\ldots$

Following analogous arguments as in the proof of Lemma \ref{lemm.H}, we can see that , $\forall\; 1\leq i\leq 2k-2$, $c_{i1}=d_{i1}$. For $i=2k-1$, we get:
$$c_{2k-1,1}=a_{2k-1,1}+\sum_{j=k}^{2k-1}a_{2k-1,j}w_j$$

\noindent Now using Theorem \ref{prop.gknested} and that $a_{2k-1,1}=v_{2k-1}$, this yields:
$$c_{2k-1,1}=v_{2k-1}+w_{2k-1}+kv_kw_k+\sum_{j=k+1}^{2k-2} (2k-j-1)v_{2k-j}w_j $$

\noindent Finally, since $(1,v)\in \mathcal G_k$,  $c_{2k-1,1}=v_{2k-1}+w_{2k-1}+kv_{k}w_{k}$.





\hspace{15cm}\end{proo}

Theorem \ref{theo.partialcommute} implies that $[\mathcal G_k,\mathcal G_k]<\mathcal G_{2k}$. As a corollary we finally obtain:

\begin{coro} \label{coro.new} For all $n\geq 1$, $\mathcal A^{(n)}\subseteq \mathcal G_{2^n}$.
\end{coro}

\begin{proo} We proceed by induction. The result is true for $n=1$: just see that if $A,B$ are lower triangular finite or infinite matrices, the elements in the main diagonal of the product $AB$, are the product of the corresponding elements in the main diagonals of $A,B$. Also the elements in the main diagonal of  $A^{-1}$ are the multiplicative inverses of the elements in the main diagonal of $A$.  So $[A,B]\in \mathcal G_2$.

 Now let us go for the case $n\geq 2$, assuming the result is true for $n-1$. If $A,B\in \mathcal A^{(n-1)}$, we have that $A,B\in \mathcal G_{2^{n-1}}$. According to the previous lemma $\Pi_{2^n-1}(A),\; \Pi_{2^n-1}(B)$ commute,  so $\Pi_{2^n-1}(A^{-1}B^{-1}AB)$ is the identity matrix, and $[A,B]\in \mathcal G_{2^n}$ by definition.

\hspace{15cm}\end{proo}

We stop now to put the above results in context with some previous works. Babenko in  \cite{BB.S}  deeply studied the \emph{substitution group} $\mathcal J$ of formal power series. He also defined, for $k\geq 2$, the groups $\mathcal J^k=(x+x^k\mathbb K[[x]])$.  In that article it is already proved, by a direct computation, that $[\mathcal J^{k},\mathcal J^k]\subset \mathcal J^{2k}$,
 and then that the $(n-1)$-th commutator $\mathcal J^{(n-1)}=[\mathcal J^{(n-2)},\mathcal J^{(n-2)}]$ is contained in $ \mathcal J^{2n}$  (Lemma 2.3, and Section 2.3 in \cite{BB.S}). 

Via the isomorphism described in Remark \ref{rem.equiv} between $\mathcal A$ and $\mathcal F_1$, we can see that the groups $\mathcal J^k$ and $\mathcal G_k$ are isomorphic. Those related results appearing in \cite{BB.S} and our Corollary \ref{coro.new} are, in some sense, equivalent. However, they have been  proved with different methods. But using Riordan matrices in the associated subgroup (Theorem \ref{main.1}) we strength the result contained in Babenko \cite{BB.S} about the substitution group, proving actually that    $\mathcal J^{(n-1)}=\mathcal J^{2n}$.



To conclude this section we prove one more result concerning the so called \emph{weighted Schr\"oder equation}. It is important for the rest of the paper and it is closely related to the groups $\mathcal G_k$.  The name of this equation comes from the analysis of the \emph{weighted composition operators}. It is the  functional equation $d u(h)=\lambda u $. In this case, it is considered in the formal power series context in the indeterminate $u=u_0+u_1t+u_2t^2+\ldots\in\mathbb K[[t]]$ and  for some given  $d\in\mathcal F_0,\;h\in\mathcal F_1$, $\lambda\in\mathbb K$.  Weighted Schr\"oder equations are suitable to be adressed in terms of Riordan matrices, because they are  eigenvector problems:
$$d u(h)=\lambda u \Longleftrightarrow (d,h) \begin{bmatrix} u_0\\ u_1 \\ u_2\\ \vdots \end{bmatrix}=\lambda  \begin{bmatrix} u_0\\ u_1 \\ u_2\\ \vdots \end{bmatrix}$$

\noindent If the previous weighted Schr\"oder equation has a solution, then $\lambda=d(0)$.  So just dividing both sides of the equation by $\lambda$, we can reduce our study, changing suitably $d$, to equations of the type:
\begin{equation}\label{eq.simp} d u(h)=u,\qquad\qquad \text{with }d(0)=1 \end{equation}

 The following Proposition  is also a detailed example of how the rest of the results in this article are proved using \emph{induction in the size of the matrices in Riordan groups of finite matrices}.

\begin{prop} \label{FE.MCWS1} $ $
\begin{itemize}
\item[(1)] Suppose $h=rx+h_2x^2+\ldots$ where $r\in\K$ and it is not a root of unity. Then there is a unique solution $u\in (1+x\mathbb K[[x]])$ of \eqref{eq.simp}, for any $d$.
\item[(2)]  Suppose $k\geq 2$ and $h=x+h_kx^k+\ldots $, with $h_k\neq 0$. Then, there exists a solution $u\in (1+x\mathbb K[[x]])$ of  \eqref{eq.simp}  if and only if $d\in (1+x^k\K[[x]])$. When this $u$ exists, it is unique.

\item[(3)] Moreover, in the previous case, for $i\geq 1$,  $u\in(1+x^i\K[[x]])$ if and only if $d\in (1+x^{i+k-1}\K[[x]])$.
\end{itemize}

\end{prop}

\begin{proo}  Let $(d,h)=(d_{ij})_{0\leq i,j<\infty}$. We have that:

\begin{itemize}

\item[(1)] Equation \eqref{eq.simp} in matricial form is:
$$(d,h)\begin{bmatrix} u_0\\ u_1 \\ \vdots \\ u_n\\ \vdots \end{bmatrix}=\begin{bmatrix}1 \\ d_{10} & r \\ \vdots & \vdots & \ddots \\ d_{n0} & d_{n1}&\ldots & r^n \\ \vdots  & \vdots &  & \vdots &\ddots\end{bmatrix}\begin{bmatrix} u_0\\ u_1 \\ \vdots \\ u_n\\ \vdots\end{bmatrix}= \begin{bmatrix} u_0\\ u_1 \\ \vdots\\ u_n\\ \vdots \end{bmatrix}$$


Now we solve the equations corresponding to each row by the method of forward substitution:
$$\begin{cases} u_0=u_0\\ d_{10} u_{0}+ru_1=u_1\\
\vdots \\ d_{n0}u_0+d_{n1}u_1+\ldots+d_{n-1,n}u_{n-1}+r^nu_n=u_n\\ \vdots \end{cases} $$

From the fact that  $r^n\neq 1$ for all $n$, we have that the infinite system above has  solution. In fact, for each  $u_0$ given there is a unique solution of the system. Choose $u_0=1$.

\item[(2)]  If there is a solution, then necessarily $d\in(1+x^k\K[[x]])$ because if $0<i<k$, then:
$$\begin{bmatrix}1 \\ 0 & 1 \\ \vdots & \ddots& \ddots \\ 0 & \hdots & 0 & 1  \\ d_{i0} &0 & \hdots & 0 & 1\end{bmatrix}\begin{bmatrix}1 \\ u_1 \\ \vdots \\ u_{i-1}\\ u_i \end{bmatrix}=\begin{bmatrix}1 \\ u_1 \\ \vdots \\ u_{i-1} \\ u_i\end{bmatrix}\Rightarrow d_{i0}=0 $$

To prove the converse, we show by induction over $n$, using the projections $\Pi_n(d,h)=(d,h)_n$, that there exist a unique sequence $(1,u_1,\ldots,u_{n-k+1})$ such that:
\begin{equation} \label{SE.EZ} (d,h)_n\begin{bmatrix}1 \\ u_1 \\ \vdots \\ u_n \end{bmatrix}=\begin{bmatrix}1 \\ u_1 \\ \vdots \\ u_n \end{bmatrix} \end{equation}

\noindent independently on the values of the terms in the sequence  $(u_{n-k+2},\ldots, u_{n})$. According to Theorem \ref{prop.gknested}, if we write \eqref{SE.EZ} as a system of linear equations in the indeterminates $u_1,\ldots,u_n$, then the values $u_{n-k+2},\ldots, u_{n}$ do not appear in them.

The first case, $n= k-1$, is obviously true since $(d,h)_{n}=I_{n+1}$ (identity matrix). Now consider the case  $n>k-1$. By induction hypothesis  there exists a unique sequence $(1,u_1,\ldots,u_{n-k})$ such that:
\begin{equation} \label{SE.EZn-1} (d,h)_{n-1}\begin{bmatrix}1 \\ u_1 \\ \vdots \\ u_{n-1} \end{bmatrix}=\begin{bmatrix}1 \\ u_1 \\ \vdots \\ u_{n-1} \end{bmatrix}\end{equation}

\noindent independently on the sequence  $(u_{n-k+1},\ldots, u_{n-1})$. The system \eqref{SE.EZ} has only one more linear equation than \eqref{SE.EZn-1}. So, we only need to check that there exists a unique $u_{n-k+1}$ such that the last equation  in \eqref{SE.EZ} holds:
$$d_{n0}+\left(\sum_{m=1}^{n-k}d_{nm}u_m\right)+d_{n,n-k+1} u_{n-k+1}+u_n=u_n$$

\noindent The elements  $d_{ij}$ with $i-j=k-1$ form an arithmetic progression, according to Theorem \ref{prop.gknested}. Recall that the first term is 0 and its common difference is $d_{k1}=h_k\neq 0$. So, we can see that this equation is equivalent to:
\begin{equation}\label{SE.CZ}(n-k+1)h_ku_{n-k+1}=-d_{n0}-\left(\sum_{m=1}^{n-k}d_{nm}u_m\right)\end{equation}

 \noindent  Then, there exists a unique value $u_{n-k+1}$ such that this equation holds. See that the right hand side in the equation above is  completely determined by the previous step.



\item[(3)] The case $i=1$ have already been discussed in (2). Now to see the general case $i>1$ see that in \eqref{SE.CZ}, taking $i=n-k+1$, 
$$ih_ku_{i}=-d_{i+k-1,0}-\left(\sum_{m=1}^{i-1}d_{i+k-1,m}u_m\right)$$

By hypothesis,  $u\in(1+x^{i-1}\mathbb K[[x]])$ and $d\in(1+x^{i+k-2}\mathbb K[[x]])$, so $u_i\neq 0$ if and only if $d_{i+k-1,0}\neq 0$.




\end{itemize}

\hspace{15cm}\end{proo}

\section{The Derived Series of  the associated subgroup of the Riordan group} \label{sec.fps}

The main result in this section is the following:

 \begin{theo}\label{main.1bis} Let $\mathcal A$ be the associated  subgroup of the Riordan group $\mathcal R$. Then for $n\geq 1$:
\begin{equation} \label{equation.main.1} \mathcal A^{(n)}=\mathcal G_{2^n}\end{equation}

\noindent  and all of its elements are commutators of elements in $\mathcal A^{(n-1)}$.

\end{theo}

Via the isomorphism between $\mathcal A$ and $\mathcal F_1$ explained in Remark \ref{rem.equiv}, this result implies Theorem \ref{main.1} below.

The case $n=1$ in the theorem states that $\mathcal A'=\mathcal G_2$, that is, $\mathcal A'$ is the subgroup of matrices of $\mathcal A$ with all the entries in the main diagonal equal to one. We will start the proof of Theorem \ref{main.1bis} showing this fact. The statement of the analogous result for formal power series is  that $\mathcal F_1'=\mathcal J$.

 \begin{lemm} \label{lemm.primer} $\mathcal A'=\mathcal G_2$.  Moreover, any element in $\mathcal G_2$ is a commutator of elements in  $\mathcal A$.

\end{lemm}

 \begin{proo}  In Corollary \ref{coro.new} we already proved that $\mathcal A'\subset \mathcal G_2$. So, the only thing we need to see is that any element $(1,g)\in\mathcal G_2$ (that is, with $g=x+g_2x^2+g_3x^3+\ldots$ or equivalently an element in $\mathcal A$ with all the entries in the main diagonal equal  to 1) is a commutator of elements in $\mathcal A$.

We will prove a stronger fact. Fix this $(1,g)$. For any $ r\in\K$ not a root of unity, there exists a unique $(1,v)\in\mathcal G_2\subset \mathcal A$ (and so with all the entries in its main diagonal equal to 1)  such that:
$$(1,g)=(1,rx)^{-1}(1,v)^{-1}(1,rx)(1,v) $$

\noindent Note that $(1,rx)$ is a diagonal Riordan matrix.


Using the inverse limit approach to the Riordan group (see again \cite{LMMPS.IL}), this is equivalent to showing that there exists a unique  sequence of finite Riordan matrices $(B_0,B_1,B_2, \ldots,B_i,\ldots)$ satisfying the following conditions. 

For all $i\geq 0$, $B_i\in\Pi_i(\mathcal A)$, $B_i=p_{i}(B_{i+1})$, $B_i$ has all the entries in the main diagonal equal to one and $(1,g)_i=(1,rx)_i^{-1}B_i^{-1}(1,rx)_iB_i $.

 Recall that  any element in $\Pi_i(\mathcal A)$ is determined by its second column. So, if we denote by $B_i=(b_{lm})_{0\leq l,m< i}$, to solve the last equation in the previous paragraph is equivalent to solve:
\begin{equation} \label{eq.lemmprimer} B_i (1,rx)_i\begin{bmatrix} 0 \\ 1 \\ g_2 \\ \vdots \\ g_i \end{bmatrix}=(1,rx)_i\begin{bmatrix} 0 \\ 1 \\ b_{21} \\ \vdots \\b_{i1}  \end{bmatrix} \end{equation}

For $i=1$, \eqref{eq.lemmprimer} is a trivial equality, independently on $r$.  The first case with special meaning is $i=2$. We have:
$$\begin{bmatrix} 1 \\ 0 & 1 \\ 0 & b_{21} & 1\end{bmatrix} \begin{bmatrix} 1 \\ 0 & r \\ 0 & 0 & r^2\end{bmatrix} \begin{bmatrix}  0 \\ 1 \\ g_2 \end{bmatrix}=\begin{bmatrix} 1 \\ 0 & r \\ 0 & 0 & r^2\end{bmatrix}\begin{bmatrix} 0 \\ 1 \\ b_{21}\end{bmatrix} $$

\noindent This matricial equation is equivalent to a system of three linear equations. The first two of them correspond to the case $i=1$. The last one is:
$$rb_{21}+r^2g_2=r^2b_{21}$$

\noindent and, since $r\neq 1$, then the solution is:
$$b_{21}=\frac{r}{r-1}g_{2}. $$

In the general case $i>2$, we need to show that there is a unique  $b_{i1}$ satisfying the last equation in the system \eqref{eq.lemmprimer}. Note that the previous case have already fixed the unique possible choice of $b_{21},\ldots,b_{i-1,1}$ satisfying the rest of the equations. This last equation is:
$$\begin{bmatrix} 0 & b_{i1} & \ldots & b_{ii} \end{bmatrix}(1,rx)_i\begin{bmatrix}0 \\ 1 \\ g_2\\ \vdots \\ g_i \end{bmatrix}=\begin{bmatrix} 0 & \ldots & 0 & r^i\end{bmatrix}\begin{bmatrix} 0 \\ 1 \\ b_{21} \\ \vdots \\ b_{i1} \end{bmatrix} \Longrightarrow$$
$$\Longrightarrow rb_{i1}+\sum_{k=2}^i\left(b_{ik} r^kg_k\right)=r^i b_{i1} \Longrightarrow (1-r^{i-1})b_{i1}=-\sum_{k=2}^i\left(b_{ik} r^{k-1}g_k\right) $$

We have to recall at this point  that, for $2\leq m\leq i$,  all the entries $b_{im}$ are determined by the elements $b_{l1}$ with $1\leq l\leq i-1$, because $B_i$ is a finite Riordan matrix. Thus, they are determined by induction hypothesis. So, nothing in  the right hand side of the last equation depends on $b_{i1}$ and all the terms are known. Since the coefficient $(1-r^{i-1})$  is not 0 ($r$ is not a root of unity) the equation has a unique solution, and the proof is complete.

\hspace{15cm}\end{proo}




Finally, we are ready to prove the main result of this section:
\medskip

\begin{proo}[Proof of Theorem \ref{main.1bis}] We have already proved that $\mathcal A'=\mathcal G_2$. Let us proceed by induction. Assume that $\mathcal A^{(n-1)}=\mathcal G_{2^{n-1}}$.  Corollary \ref{coro.new} guarantees that  $\mathcal A^{(n)}\subseteq \mathcal G_{2^{n}}$. So, we only need to prove that every element in $\mathcal G_{2^n}$ is a commutator  of elements in $\mathcal A^{(n-1)}$.

Again, we will prove a stronger fact. Let us fix an element $(1,g)\in\mathcal G_{2^n}$ and suppose that $A=(1,x+\lambda x^{2^{n-1}})$ with $\lambda\neq 0$. In particular, $A$ in $\mathcal A^{(n-1)}$ and not in $\mathcal G_{2^{n-1}+1}$. Then, we have to prove that there exists a unique $B\in \mathcal G_{2^{n-1}+1}$ such that:
$$(1,g)=A^{-1}B^{-1}AB $$

As in the previous lemma, this is equivalent to prove that

\noindent \textbf{Claim (n):} Fixed $n$, there exists a unique sequence of finite Riordan matrices $(B_0,B_1,B_2,\ldots,B_i,\ldots)$ satisfying the following conditions. $B_{2^{n-1}}$ must be the identity matrix of the corresponding size and for all $i\geq 0$, $B_i\in\Pi_i(\mathcal A)$, $B_i=p_i(B_{i+1})$. Finally, for each $i$, $B_i$ must satisfy:
\begin{equation} \label{eq.theoasociado} B_i\Pi_i(A) (1,g)_i=\Pi_i(A) B_i\end{equation}





 Denote $A=(a_{lm})_{0\leq l,m<\infty}$, $B_i=(b_{lm})_{0\leq l,m<i}$. For $i\geq 2^n$, Equation \eqref{eq.theoasociado} is equivalent to
$$\left[\begin{array}{c   | c}  &  \\  B_{i-1}  & \\ & \\ \hline b_{i0} \hdots b_{i,i-1} & 1 \end{array} \right]\left[\begin{array}{c   | c}   & \\ \Pi_{i-1}(A)  & \\  & \\ \hline a_{i0} \hdots a_{i,i-1} & 1 \end{array} \right]    \begin{bmatrix} 0\\ 1 \\ 0 \\ \vdots \\ 0 \\ g_{2^n} \\ \vdots \\ g_i  \end{bmatrix} =\left[\begin{array}{c  | c}  &  \\  \Pi_{i-1}(A)  & \\  & \\ \hline a_{i0} \hdots a_{i,i-1} & 1 \end{array} \right] \begin{bmatrix} 0 \\ 1 \\ 0 \\ \vdots \\ 0 \\ b_{2^{n-1},1}\\ \vdots \\ b_{i1}\end{bmatrix}   $$

Note that the last equation in the above linear system, in the indeterminates $b_{2^{n-1},1},\ldots,b_{i1}$, is:
\begin{equation}\label{DS.GE}\underbrace{[b_{i0} \hdots b_{i,i-1}] \Pi_{i-1}(A)\begin{bmatrix}0 \\ 1 \\ 0 \\ \vdots \\ 0 \\ g_{2^n} \\ \vdots \\ g_{i-1} \end{bmatrix}}_{(I)}+[ a_{i0}\hdots a_{i,i-1} ]\begin{bmatrix}0 \\ 1 \\ 0 \\ \vdots \\ 0 \\ g_{2^n} \\ \vdots \\ g_{i-1} \end{bmatrix}+g_i =\underbrace{[a_{i0} \hdots a_{i,i-1} ]\begin{bmatrix} 0 \\ 1 \\ 0 \\ \vdots \\ 0 \\ b_{2^{n-1},1}\\ \vdots \\ b_{i-1,1}\end{bmatrix}}_{(II)} +b_{i1} \end{equation}


\noindent Appart from $b_{i1}$, only the terms (I) and (II) in the above equation contain indeterminates.


Now we will see, as happened in some of the previous proofs and because of the structure of the matrices, that in the  linear system above only the unkwnons $b_{11},\ldots,b_{i-2^{n-1}+1,1}$ are involved.

First,  note that $A=(1,x+\lambda x^{2^{n-1}})\in \mathcal G_{2^{n-1}}$. So, according to Theorem \ref{prop.gknested}, we obtain
\begin{equation}  \label{eq.t7par} a_{j+m,j}=\begin{cases}1&\text{if } m=0\\ 0 &\text{if } 1\leq m\leq 2^{n-1}-2\\
j \lambda &\text{if }m=2^{n-1}-1\\
 0 &\text{if } 2^{n-1}\leq m\leq 2^n-2 \end{cases} ,\qquad\qquad b_{j+m,j}=\begin{cases}1&\text{if } m=0\\ 0 &\text{if } 1\leq m\leq 2^{n-1}-1\\
 j b_{j+1,1} &\text{if } 2^{n-1}\leq m\leq 2^n-2 \end{cases}\end{equation}

We have that:
$$(I)=\left[0,\; b_{i1}+ b_{i,2^{n-1}}\lambda,\; b_{i2}+\sum_{k=2^{n-1}+1}^{i-1} b_{ik}a_{k2},\; \ldots\;,\; b_{i,i-1}+\sum_{k=2^{n-1}+i-2}^{i-1} b_{ik}a_{k,i-1}\right]\begin{bmatrix}0 \\ 1 \\ 0 \\ \vdots \\ 0 \\ g_{2^n} \\ \vdots \\ g_{i-1} \end{bmatrix}$$

\noindent Note that all the elements $b_{ik}$ in any of the sums $\sum_{k=2^{n-1}+1}^{i-1} b_{ik}a_{k2}$, \ldots, $\sum_{k=2^{n-1}+i-2}^{i-1} b_{ik}a_{k,i-1}$ depend only on $b_{11},\ldots, b_{i-2^{n-1},1}$. So:
\begin{equation}\label{eq.C1}(I)=b_{i1}+\lambda b_{i,2^{n-1}}+[C_1]\end{equation}

\noindent  where nothing in $[C_1]$ depends on $b_{i-2^{n-1}+1,1}, b_{i-2^{n-1}+2,1},\ldots$

Moreover,
\begin{equation} \label{eq.C} (II)=[a_{i0} \hdots a_{i,i-1} ]\begin{bmatrix} 0 \\ 1 \\ 0 \\ \vdots \\ 0 \\ b_{2^{n-1},1}\\ \vdots \\ b_{i-1,1}\end{bmatrix} =\left[\sum_{k=1}^{i-2^{n-1}}a_{ik}b_{k1}\right]+ a_{i,i-2^{n-1}+1}b_{i-2^{n-1}+1,1}\end{equation}

 \noindent where, again, nothing inside the term $\left[\sum_{k=1}^{i-2^{n-1}}a_{ik}b_{k1}\right]$  depends on $b_{i-2^{n-1}+1,1}, b_{i-2^{n-1}+2,1},\ldots$

\noindent Substituting \eqref{eq.C1} and \eqref{eq.C} in \eqref{DS.GE}, reorganizing and cancelling when needed, we obtain
\begin{equation}\label{DS.GE2}  a_{2^{n-1},1}b_{i,2^{n-1}}-a_{i,i-2^{n-1}+1}b_{i-2^{n-1}+1,1}=\underbrace{-[C_1]-[ a_{i0}\hdots a_{i,i-1} ]\begin{bmatrix}0 \\ 1 \\ 0 \\ \vdots \\ 0 \\ g_{2^n} \\ \vdots \\ g_{i-1} \end{bmatrix}-g_i +\left[\sum_{k=1}^{i-2^{n-1}}a_{ik}b_{k1}\right]}_{[C_2]}\end{equation}

\noindent The right hand side in \eqref{DS.GE2}, i.e. $[C_2]$, does not depend on $b_{i-2^{n-1}+1,1}$, $b_{i-2^{n-1}+2,1},\ldots$ 

Following \eqref{eq.t7par} we get $a_{2^{n-1},1}=\lambda$, $a_{i,i-2^{n-1}+1}=(i-2^{n-1}+1) \lambda$. On the other hand, by definition of finite Riordan matrix, there exist a matrix $B=(1,h)\in\mathcal A$ such that $\Pi_i(B)=B_i$. The entries in the columns in $B$ are the coefficients of a geometric sequence with first term the formal power series $1$ with common ratio the formal power series $h$. Since $b_{i-2^{n-1}+1,1}=[x^{i-2^{n-1}+1}]h$ and $ b_{i,2^{n-1}}=[x^i]h^{2^{n-1}}$, then
$$b_{i,2^{n-1}}=b_{i-2^{n-1}+1,1}+[C_3]$$

\noindent  where nothing in the term $[C_3]$ depends on $b_{i-2^{n-1}+1,1}, b_{i-2^{n-1}+2,1},\ldots$

So, finally, \eqref{DS.GE2} is equivalent to:
\begin{equation} \label{eq.finalring}  \lambda  (b_{i-2^{n-1}+1,1}+[C_3])-(i-2^{n-1}+1) \lambda b_{i-2^{n-1}+1,1}=[C_2]\end{equation}

\noindent which has a unique solution in the indeterminate $b_{i-2^{n-1}+1,1}$, provided that $b_{k1}$ is known for $1\leq k\leq i-2^{n-1}$.

We can now prove  Claim (n) by induction over $i$. Let us begin imposing that $B_{2^{n-1}}$ is the identity matrix. We trivially have (Theorem \ref{theo.partialcommute}) that \eqref{eq.theoasociado} holds for $i\leq 2^n-1$. Assume that $B_{i-1} \Pi_{i-1}(A)(1,g)_{i-1}=\Pi_{i-1}(A) B_{i-1}$ for some $i\geq 2^n$ for a fixed value of the intedeterminates $b_{21},\ldots,b_{i-2^{n-1},1}$ indendently on the value of the indeterminates $b_{i-2^{n-1}+1,1}, b_{i-2^{n-1}+2,1},\ldots$ Then, as explained before, Equation \eqref{eq.theoasociado} holds if and only if its last equation holds. We have just showed that there is a unique solution $b_{i-2^{n-1}+1,1}$ of this equation and that this solution does not depend on $b_{i-2^{n-1}+2,1}, b_{i-2^{n-1}+3,1},\ldots$

\hspace{15cm}\end{proo}

 \begin{theo}\label{main.1} For $n\geq 1$:
\begin{equation} \label{equation.main.1}  \mathcal F_1^{(n)}= \mathcal J^{(n-1)}=\{g\in \mathcal J: g\in(x+x^{2^{n}}\mathbb K[[x]])\}\end{equation}

\noindent  and all of its elements are commutators in $\mathcal F_1^{(n-1)}$.

\end{theo}

\begin{rem} \label{rem.move} In the case in which $\mathbb K$ is a finite field, there is an analogue of Theorem \ref{main.1}, in the context of the so called Nottingham groups (see \cite{C.NG}). Consequently, we could also obtain an analogue of Theorem \ref{main.1bis}.

\end{rem}

\section{Proof of Theorem \ref{main.2} and some consequences} \label{sec.WSE}

Finally we can prove Theorem \ref{main.2}  as a  consequence of Proposition \ref{FE.MCWS1} and Theorem \ref{main.1}.

\medskip

\begin{proo}[Proof of Theorem \ref{main.2}] We are going to prove by  induction over $n$ that:
$$ \mathcal R^{(n)}=\begin{cases}\mathcal R & \text{if }n=0\\
\{(d,h)\in\mathcal R: d\in(1+x^{2^n-n}\mathbb K[[x]]),h\in (x+x^{2^n}\mathbb K[[x]])\} & \text{if }n>0\end{cases}$$

The case $n=0$ is just a notational fact. Let us prove the case $n>0$, assuming that the case $n-1$ holds. The set:
$$\{(d,h)\in\mathcal R: d\in(1+x^{2^n-n}\mathbb K[[x]]),h\in (x+x^{2^n}\mathbb K[[x]])\}$$

\noindent is a group. So we only need to prove that any commutator of elements in $\mathcal R^{(n-1)}$ is of this type.

Using the multiplication and inversion formula in the Riordan group, se have  that the matricial equation $ (d,h)=(u,v)^{-1}(f,g)^{-1}(u,v)(f,g)$  is equivalent to the following equation in formal power series
\begin{equation} \label{eq.condition} \begin{cases}h=g\circ v\circ g^{-1}\circ v^{-1} \\ \left[\frac{u}{u(g) d(v(g))}\right] f(v)=f \end{cases}\end{equation}

If $v,g\in(x+x^{2^{n-1}}\mathbb K[[x]])$ we have already proved that $h\in(x+x^{2^n}\mathbb K[[x]])$ (consequence of Theorem \ref{theo.partialcommute}), and that in this case, for any fixed $h$ we can find such a pair $v,g$ (Theorem \ref{main.1}).

On the other hand, according to Proposition \ref{FE.MCWS1}, there exists $f\in(1+x^{2^{n-1}-n+1}\mathbb K[[x]])$, satisfying the second condition in \eqref{eq.condition} if and only if
\begin{equation} \label{eq.prodcoeff}\left[\frac{u}{u(g)d(v(g))}\right]\in (x+x^{2^n-n}\mathbb K[[x]])\end{equation}

\noindent  To check \eqref{eq.prodcoeff}, see that, as soon as $u\in(1+x^{2^{n-1}-n+1}\mathbb K[[x]]) $ we obtain that $\frac{u}{u(g)}\in (1+x^{2^n-n}\mathbb K[[x]])$. Note that $u$ and $u(g)$ only differ from the $(x^{2^n-n})$-th term onwards.  On the other hand, $\frac{1}{ d(v(g))}\in(x+x^{2^n-n}\mathbb K[[x]])$ if and only if $d\in(x+x^{2^n-n}\mathbb K[[x]])$. As a consequence,   \eqref{eq.prodcoeff} holds if and only if $d\in(x+x^{2^n-n}\mathbb K[[x]])$ and the proof is finished.

\hspace{15cm}\end{proo}

\noindent The case $n=1$ in Theorem \ref{main.2} was previously proved in \cite{LMP.GGI}.

Recall that $\Pi_k(\mathcal R)=\mathcal R_k$ and $\Pi_k$ is a group homomorphism. So, using Lema 2.1 in \cite{N} we get

\begin{coro}
   For every $n,k \geq 1$, $\mathcal R_k^{(n)}=\Pi_k(\mathcal R^{(n)})$. 
\end{coro}

\begin{rem}

\begin{itemize}
\item[(a)] Note that, as soon as $k<2^n-n$, $\mathcal R_k^{(n)}$ is formed only by the identity matrix in $\mathcal R_k$.

\item[(b)]  If $2^n-n\leq k < 2^n$,  $\mathcal R_k^{(n)}$ is formed only by Toeplitz matrices and then in this case $\mathcal R_k^{(n)}$ is abelian.

\item[(c)] Note also that $\mathcal R_k$ is  solvable.

\end{itemize}

\end{rem}

\begin{rem}\textbf{\emph{(Some prospects)}} For any integer $k\geq 0$ we denote by $\varphi(k)$ to the \emph{derived length} or \emph{solvable length} of $\mathcal R_k$, that is, the lowest integer $n$ such that $\mathcal R_k^{(n)}$ is the trivial subgroup.  The first terms of the sequence $\varphi$ are
$$1, 2, 3, 3, 3, 4, 4, 4, 4, 4, 4, 4, 5,\ldots $$

\noindent They are the first terms of the sequence A103586  in  OEIS. See \cite{D, GLASBY, NBF}, and the references therein, for some featured properties of the solvable length in special cases.



A description of the derived series of the Riordan group has not been given yet for other choices of the field $\mathbb K$ (cases in which $\mathbb K$ is not of characteristic 0). See  Remark \ref{rem.move}.

For a better understanding of the algebraic structure of the Riordan group it would be desirable to determine the quotient groups $\mathcal R^{(n)}/\mathcal R^{(n+1)}$. It is easy to see, and it is proved in \cite{LMP.GGI}, that $\mathcal R/\mathcal R'$ is isomorphic to $\mathbb K^*\times\mathbb K^*$. But the rest of the elements in this sequence of quotients do not seem to have such an easy description.

We would like to point out the following. Let $\mathfrak{R}$ be a commutative ring with a unit, denoted by 1. Consider a coherent definition of Riordan matrix with entries in $\mathfrak{R}$ and such that all  the entries in the main diagonal equal to $1$. The set of such matrices should be also a group (with a certain operation). For this group, we can study the corresponding derived subgroup. The case where $\mathfrak{R}=\mathbb Z$ (the ring of integers) is of special interest in combinatorics. It would be interesting, due to the possible applications, to compute the derived series. The arguments followed herein do not apply for rings, since we usually need to take multiplicative inverses. See \eqref{eq.finalring} and \eqref{SE.CZ}.

\end{rem}


\begin{thebibliography}{9}













\bibitem{BB.S}\textit{I. K. Babenko}, {Algebra, geometry, and topology of the substitution group of formal power series}, Russian
Math. Surveys \textbf{68} (2013), 1--68.




\bibitem{BBN} \textit{M. Barnabei, A. Brini} and \textit{G. Nicoletti}, Recursive matrices and umbral calculus, J. Algebra \textbf{75} (1982), 546--573.

\bibitem{BH} \textit{P. Barry} and \textit{A. Hennessy}, Meixner-type results for Riordan arrays and associated integer sequences, J. Integer Seq. \textbf{13} (2010).



\bibitem{C.NG} \textit{R. Camina} {Notingham Group} in {New Horizons in pro-p Groups}, Progress in Mathematics, vol. 184, Springer Science and Business Media, Boston (2012).

\bibitem{C} \textit{P. J. Cassidy} Products of commutators are not always commutators: an example, Amer. Math. Monthly \textbf{86} (1979), 772--772.


\bibitem{CLW} \textit{X. Chen, H. Liang} and \textit{Y. Wang}, Total positivity of Riordan arrays, European J. Combin. \textbf{46} (2015), 68--74.






\bibitem{C.G1} \textit{G. S. Cheon}, \textit{J. H. Jung}, \textit{S. Kitaev} and \textit{S. A. Mojallal} Riordan graphs I: Structural properties, Linear Algebra Appl. \textbf{579} (2019) 89--135.






\bibitem{CLMPS.L} \textit{G. S. Cheon}, \textit{A. Luz\'on}, \textit{M. A. Mor\'on}, \textit{L. F. Prieto-Mart\'inez} and \textit{M. H. Song}, {Finite and infinite dimensional Lie group structures on Riordan groups}, Adv. Math. \textbf{319} (2017), 522--566.

\bibitem{D} \textit{J. D. Dixon}, The solvable length of a solvable linear group, Math. Z. \textbf{107} (1968), 151--158.

\bibitem{GLASBY} \textit{S. P. Glasby}, The composition and derived lengths of a soluble group, J. Algebra \textbf{120} (1989), 406--413.


\bibitem{GO} \textit{C. K. Gubta} and \textit{W. Holubowski}, Commutator subgroup of Vershik-Kerov group, Linear Algebra Appl. \textbf{436} (2012), 4279--4284.

\bibitem{G} \textit{R. M. Guralnick}, Commutators and commutator subgroups, Adv. Math. \textbf{45} (1982), 319--330.


\bibitem{HHS} \textit{T. X. He}, \textit{L. C. Hsu} and \textit{P. J. S. Shiue}, The Sheffer group and the Riordan group, Discrete Appl. Math. \textbf{155} (2007), 1895--1909.


\bibitem{H} \textit{I. C. Huang}, Inverse relations and Schauder bases, J. Combin. Theory, Ser. A \textbf{97} (2002), 203--224.

\bibitem{I} \textit{I. M. Isaacs}, Commutators and the commutator subgroup, Amer. Math. Monthly \textbf{84}, 720--722.

\bibitem{JAB1} \textit{E. Jabotinsky}, Sur la representation de la composition de fonctions par un produit de matrices, Application d l'iteration de $e^t$ et de $e^t-1$, C. R. Acad. Sci. Paris \textbf{224} (1947), 323--324.


\bibitem{JAB2} \textit{E. Jabotinsky}, Representation of functions by matrices. Application to Faber polynomials, Proc. Amer. Math. Soc.  \textbf{4} (1953), 546--553.


\bibitem{JN} \textit{C. Jean-Louis} and \textit{A. Nkwanta}, Some algebraic structure of the Riordan group, Linear Algebra Appl. \textbf{438} (2013), 2018--2035.

\bibitem{JENN} \textit{S. A. Jennings}, {Substitution groups of formal power series}, Canad. J. Math.  \textbf{6} (1954), 325--340.










\bibitem{LIEBECK} \textit{M. W. Liebeck}, \textit{E. A. O'Brien}, \textit{A. Shalev} and \textit{P. H. Tiep}, The Ore conjecture, J. Eur. Math. Soc. \textbf{12} (2010), 939--1008.

\bibitem{L} \textit{A. Luz\'on}, {Iterative processes related to Riordan arrays: The reciprocation and the inversion of power series}, Discrete mathematics \textbf{310} (2010), 3607--3618.



\bibitem{LM.U} \textit{A. Luz\'on}, \textit{M. A. Mor\'on}, {Ultrametrics, Banach's fixed point theorem and the Riordan group}, Discrete Appl. Math. \textbf{156} (2008), 2620--2635.



\bibitem{LMMPS.IL} \textit{A. Luz\'on}, \textit{D. Merlini}, \textit{M. A. Mor\'on}, \textit{L. F. Prieto-Mart\'inez} and \textit{R. Sprugnoli}, {Some inverse limit approaches to the Riordan group}, Linear Algebra  Appl. \textbf{491} (2016), 239--262.







\bibitem{LMP.I} \textit{A. Luz\'on}, \textit{M. A. Mor\'on} and \textit{L. F. Prieto-Mart\'inez}, {A formula to construct all involutions in Riordan matrix groups}, Linear Algebra and its Applications \textbf{533} (2017), 397--417.

\bibitem{LMP.GGI} \textit{A. Luz\'on}, \textit{M. A. Mor\'on} and \textit{L. F. Prieto-Mart\'inez}, {The group generated by Riordan involutions}, Rev. Mat. Complut., \url{doi.org/10.1007/s13163-020-00382-8}.

\bibitem{LMR} \textit{A. Luz\'on, M. A. Mor\'on} and \textit{J. L. Ram\'irez}, On Ward's differential calculus, Riordan matrices and Sheffer polynomials, Linear Algebra Appl. \textbf{610}, 440--473.

\bibitem{MRSV} \textit{D. Merlini, D. G. Rogers, R. Sprugnoli} and \textit{M. C. Verri}, On some alternative characterizations of Riordan arrays, Canad. J.  Math. \textbf{49}, (1997), 301--320.

\bibitem{N} \textit{B. H. Neumann}, On ascending derived series, Compos. Math. \textbf{13}, (1956), 47--64.

\bibitem{NBF} \textit{M. F. Newman}, The Soluble Length of Soluble Linear Groups, Math. Z. \textbf{126} (1972), 59--70.

\bibitem{O.C} \textit{O. Ore}, {Some remarks on commutators}, Proc. Amer. Math. Soc. \textbf{2} (1951), 307--314.

\bibitem{P} \textit{L. F. Prieto-Mart\'inez}, {Geometric Continuity of plane curves in terms of Riordan matrices and an application to the F-chordal problem}, Rev. R. Acad. Cienc. Exactas. F\'is. Nat. Ser. A Mat. \textbf{115} (2021).


\bibitem{ROGERS} \textit{D. G. Rogers}, {Pascal triangles, Catalan numbers and renewal arrays}, Discrete Math. \textbf{22} (1978), 301--310.








\bibitem{S.C} \textit{L. W. Shapiro}, A Catalan triangle, Discrete Math. \textbf{14} (1976), 83--90.

\bibitem{Sh.OQ}\textit{ L. W. Shapiro}, {Some Open Questions about Random Walks, Involutions, Limiting Distributions, and Generating
Functions}, Adv. in Appl. Math. \textbf{27} (2001), 585--596.

\bibitem{SHAPIRO.ORIGINAL} \textit{L. Shapiro}, \textit{S. Getu}, \textit{W. J. Woan} and \textit{L. C. Woodson}, The Riordan Group, Discrete Appl. Math. \textbf{34} (1991), 229--239.





\bibitem{S} \textit{R. Slowik}, {The lower central series of subgroups of the Vershik-Kerov group}, Linear Algebra Appl.  \textbf{436} (2012), 2299--2310.








\bibitem{S.CS} \textit{R. Sprugnoli}, {Riordan arrays and combinatorial sums}, Discrete Math. \textbf{132} (1994), 267--290.

\bibitem{Th} \textit{R. C. Thompson}, {Commutators in the special and general linear groups}, Trans. Amer. Math. Soc.  \textbf{33} (1961), 16--33.

\bibitem{VS} \textit{L. Verde-Star}, {Dual operators and Lagrange inversion in several variables}, Adv. Math. \textbf{58} (1985), 89--108.

\bibitem{WW} \textit{W. Wang} and \textit{T. Wang}, Generalized Riordan arrays, Discrete Math. \textbf{308} (2008), 6466--6500.


\bibitem{Y} \textit{S. L. Yang, Y. N. Dong, L. Yang} and \textit{J. Yin}, Half of a Riordan array and restricted lattice paths, Linear Algebra Appl. \textbf{537} (2018), 1--11.

\bibitem{Z} \textit{S. Zemel}, Generalized Riordan groups and operators on polynomials, Linear Algebra Appl. \textbf{494} (2016), 286--308.

















\end{thebibliography}
\end{document}